\documentclass[12pt]{article}
\usepackage{amssymb,latexsym,pstcol,pst-plot}

\newcommand{\ben}{\begin{enumerate}}
\newcommand{\een}{\end{enumerate}}
\newcommand{\ble}{\begin{lem}}
\newcommand{\ele}{\end{lem}}
\newcommand{\bth}{\begin{thm}}
\renewcommand{\eth}{\end{thm}}
\newcommand{\bpr}{\begin{prop}}
\newcommand{\epr}{\end{prop}}
\newcommand{\bco}{\begin{cor}}
\newcommand{\eco}{\end{cor}}
\newcommand{\bcon}{\begin{conj}}
\newcommand{\econ}{\end{conj}}
\newcommand{\bde}{\begin{defn}}
\newcommand{\ede}{\end{defn}}
\newcommand{\bex}{\begin{exa}}
\newcommand{\eex}{\end{exa}}
\newcommand{\barr}{\begin{array}}
\newcommand{\earr}{\end{array}}
\newcommand{\btab}{\begin{tabular}}
\newcommand{\etab}{\end{tabular}}
\newcommand{\beq}{\begin{equation}}
\newcommand{\eeq}{\end{equation}}
\newcommand{\bea}{\begin{eqnarray*}}
\newcommand{\eea}{\end{eqnarray*}}
\newcommand{\bce}{\begin{center}}
\newcommand{\ece}{\end{center}}
\newcommand{\bpi}{\begin{picture}}
\newcommand{\epi}{\end{picture}}
\newcommand{\bpp}{\begin{picture}}
\newcommand{\epp}{\end{picture}}
\newcommand{\bfi}{\begin{figure} \begin{center}}
\newcommand{\efi}{\end{center} \end{figure}}
\newcommand{\bprf}{\begin{proof}}
\newcommand{\eprf}{\end{proof}}
\newcommand{\capt}{\caption}
\newcommand{\bsl}{\begin{slide}{}}
\newcommand{\esl}{\end{slide}}
\newcommand{\bfr}{\begin{frame}}
\newcommand{\efr}{\end{frame}}

\newcommand{\pf}{{\bf Proof}\hspace{7pt}}

\newcommand{\Qqed}{\qquad\rule{1ex}{1ex}\medskip}

\newcommand{\hso}[1]{\hspace{-1pt}}

\newcommand{\qmq}[1]{\quad\mbox{#1}\quad}

\def\<{\langle}
\def\>{\rangle}

\newcommand{\ra}{\rightarrow}

\newcommand{\bbZ}{{\mathbb Z}}

\newcommand{\cV}{{\cal V}}

\newtheorem{thm}{Theorem}[section]
\newtheorem{prop}[thm]{Proposition}
\newtheorem{cor}[thm]{Corollary}
\newtheorem{lem}[thm]{Lemma}
\newtheorem{conj}[thm]{Conjecture}
\newtheorem{exa}[thm]{Example}

\begin{document}
\pagestyle{plain}

\title{A human proof for a generalization of Shalosh B. Ekhad's $10^n$ Lattice Paths Theorem
}
\author{
Nicholas A. Loehr
\thanks{Research supported in part by an NSF Postdoctoral Research Fellowship}
\\[-5pt]
\small Department of Mathematics, College of William \& Mary\\[-5pt]
\small Williamsburg, VA, \texttt{nick@math.wm.edu}\\  
Bruce E. Sagan
\thanks{Research supported in part by DIMACS}
\\[-5pt]
\small Department of Mathematics, Michigan State University\\[-5pt]
\small East Lansing, MI, \texttt{sagan@math.msu.edu}\\
Gregory S. Warrington
\thanks{Research supported in part by an NSF Postdoctoral Research Fellowship}
\\[-5pt]
\small Department of Mathematics, Wake Forest University\\[-5pt]
\small Winston-Salem, NC, \texttt{warrings@wfu.edu}
}

\date{\today}
\maketitle

\begin{abstract}
Consider lattice paths in $\bbZ^2$ taking unit steps north ($N$) and
east ($E$).  Fix positive integers $r,s$ and put an equivalence
relation on points of $\bbZ^2$ by letting $v,w$ be equivalent if
$v-w=\ell(r,s)$ for some $\ell\in\bbZ$.  Call a lattice path {\it valid\/}
if whenever it enters a point $v$ with an $E$-step, then any further
points of the path in the class of $v$ are also entered with an
$E$-step.  Loehr and Warrington conjectured that the number of valid
paths from $(0,0)$ to $(nr,ns)$ is ${r+s\choose r}^n$.  We prove this
conjecture when $s=2$.
\end{abstract}

\section{Introduction}

A {\it lattice path\/} is a directed graph whose vertices are elements
of $\bbZ^2$ where $\bbZ$ denotes the integers.  All our lattice paths will
have edges which are unit steps north ($N$-steps) or east ($E$-steps).  It
is well-known and easy to prove that the number of such paths from
$(0,0)$ to $(r,s)$ is ${r+s\choose r}$.

Given $r,s$ we put an equivalence relation on $\bbZ^2$ by saying that
points $v,w$ are equivalent if $v-w=\ell(r,s)$ for some $\ell\in\bbZ$.
As usual, addition and scalar multiplication of points are done
componentwise.  Denote the equivalence class of $v=(x,y)$ by
$[v]=[x,y]$.  Call a 
lattice path $P$ {\it valid\/} if it satisfies the following
condition:  Whenever $P$ enters a point $v$ with an $E$-step then any
future points of $P$ in $[v]$ must be entered by an $E$-step.  
Otherwise $P$ is said to be {\it invalid}.  Figure~\ref{invalid}
shows an invalid path for $r=3$, $s=2$.  In particular, we can take
$[v]=[2,1]$
and then two points where the validity condition is violated are
shown as circles.  Loehr
and Warrington made the following conjecture.
\bcon
\label{lw}
The number of valid paths from $(0,0)$ to $(nr,ns)$ is 
$$
{r+s\choose r}^n.
$$
\econ

\bfi
\begin{pspicture}(0,0)(9,6)
\psgrid[subgriddiv=1,griddots=10](0,0)(9,6)
\psline(0,0)(1,0)
\psline(1,0)(1,1)
\psline(1,1)(2,1)
\psline(2,1)(2,4)
\psline(2,4)(8,4)
\psline(8,4)(8,6)
\psline(8,6)(9,6)
\pscircle*(2,1){.1}
\pscircle*(8,5){.1}
\rput(1.5,1.3){$E$}
\rput(8.3,4.5){$N$}
\rput(2.2,1.2){$v$}
\rput(7.8,5.2){$w$}
\end{pspicture}
\capt{An invalid path for $r=3$, $s=2$}
\label{invalid}
\efi

Ekhad, Vatter, and Zeilberger~\cite{evz:plw} gave a fully
computer-based proof of the special case $r=3$, $s=2$ of this
conjecture.  It is for this reason (and also at the request of its two
human coauthors) that we call this  ``Shalosh B. Ekhad's $10^n$ Lattice
Paths Theorem.''  Although our demonstration is very different in nature, being
purely human and bijective rather than inductive, we should mention
that some of our ideas came from looking at the trees generated by Ekhad.
Also, while we were writing this note, it came to our attention that
Jonas Sj\"ostrand~\cite{sjo:clp} has given a bijective proof of the full conjecture
which is similar to ours in some respects but differs in others.

The rest of this paper is structured as follows.  In the next section
we will provide three lemmas which will permit us to demonstrate that our
bijection is well defined.  These lemmas hold for all $r,s$.  In the
final section, we prove the case $s=2$ of Conjecture~\ref{lw}.

\section{Preliminary Lemmas}

We first need to establish some notation.  If $P$ is any lattice path
and $v$ is any point, then $P+v$ will denote the translated lattice
path obtained by adding $v$ to every point of $P$.  Note that $P$ is
invalid precisely when there is some $\ell>0$ such that $P$ and
$Q=P+\ell(r,s)$ intersect where $P$ enters the intersection
with an $N$-step and $Q$ enters it with an $E$-step.  In
Figure~\ref{invalid}, $\ell=2$.

Given $P$ and a
line $x+y=i$ we let $v_i=v_i(P)$ be the intersection of $P$ with this
line.  Note that this will not necessarily be the $i$th vertex of $P$
unless $P$ starts at the origin.  We denote the coordinates of $v_i$
by $(x_i,y_i)=(x_i(P),y_i(P))$ and the edge/step of $P$ into $v_i$ by
$e_i=e_i(P)$.

\bfi
\begin{pspicture}(0,0)(10,6)
\psgrid[subgriddiv=1,griddots=10](0,0)(10,6)
\psline(0,0)(0,1)
\psline(0,1)(6,1)
\psline(6,1)(6,4)
\psline(6,4)(7,4)
\psline[linestyle=dashed](3,2)(3,3)
\psline[linestyle=dashed](3,3)(9,3)
\psline[linestyle=dashed](9,3)(9,6)
\psline[linestyle=dashed](9,6)(10,6)
\pscircle*(4,3){.1}
\pscircle*(6,1){.1}
\pscircle*(6,3){.1}
\pscircle*(6,4){.1}
\pscircle*(7,3){.1}
\psline[arrows=<->,arrowsize=3pt 4](1,6)(7,0)
\psline[arrows=<->,arrowsize=3pt 4](4,6)(10,0)
\rput(2,6.5){$x+y=i=7$}
\rput(5,6.5){$x+y=k=10$}
\rput(.3,.7){$P$}
\rput(3.3,2.7){$Q$}
\end{pspicture}
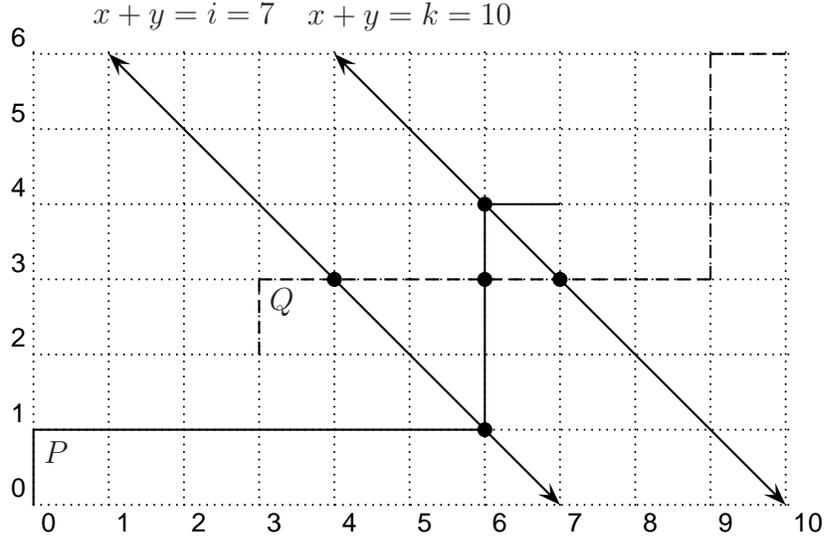
\capt{The Switching Lemma for $r=3$, $s=2$}
\label{swlfig}
\efi

The next lemma is fundamental to all that follows.  The reader may
find it useful to refer to Figure~\ref{swlfig} while reading the
statement and proof.  Circles mark the points on the lines $x+y=i$,
$j$ and $k$.
\ble[Switching Lemma]
\label{swl}
Let $P$ be a lattice path and let $Q=P+(r,s)$.  If there are integers
$i<k$ with
$$
x_i(P)>x_i(Q)\qmq{and} x_k(P)\le x_k(Q),
$$
then $P$ is invalid.
\ele
\pf\
Since the $x$-coordinate of a path changes by at most one with each
step, the hypotheses imply that there is an index $j$ with $i<j\le k$
and $x_j(P)=x_j(Q)$.  If one takes the smallest such $j$, then we must
have $e_j(P)=N$ and $e_j(Q)=E$.  It follows that $P$ is invalid by the
remark at the end of the first paragraph of this section.
\Qqed

Partially order $\bbZ^2$ componentwise, i.e., $(x,y)\le (x',y')$ if and
only if $x\le x'$ and $y\le y'$.
If $P$ is a path and $v$ is a point then we say that {\it $P$ passes
strictly west of $v$\/} if there is a point $v'$ of $P$ with the same
$y$-coordinate as $v$ satisfying $v'<v$.  
We also define $P$ to
{\it pass weakly west\/} of $v$ if for all points $v'$ of $P$
with the same $y$-coordinate we have
$v'\le v$.  (We also insist that at least one such point exists.)
Note the difference in the quantifiers between the two definitions.
Passing  east, either strictly or weakly, is defined by simply
reversing the inequalities.    Many of
our geometric arguments will be based on the following lemma.
A path satisfying the hypotheses of this result is shown in
Figure~\ref{tplfig}.

\bfi
\begin{pspicture}(0,0)(9,6)
\psgrid[subgriddiv=1,griddots=10](0,0)(9,6)
\psline(0,0)(0,1)
\psline(0,1)(3,1)
\psline(3,1)(3,2)
\psline(3,2)(6,2)
\psline(6,2)(6,5)
\psline(6,5)(8,5)
\psline(8,5)(8,6)
\psline(8,6)(9,6)
\pscircle*(2,1){.1}
\rput(2.3,.7){$u$}
\pscircle*(5,3){.1}
\rput(4.7,3.3){$v$}
\pscircle*(8,5){.1}
\rput(8.3,4.7){$w$}
\end{pspicture}
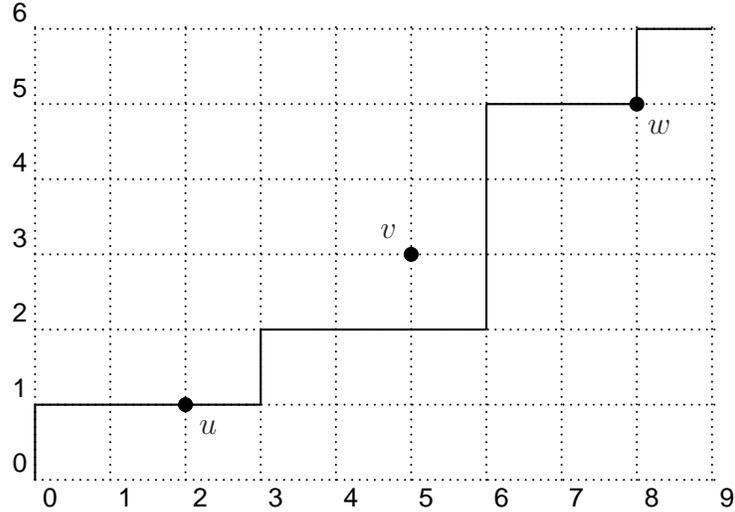
\capt{The Three-Points Lemma for $r=3$ and $s=2$}
\label{tplfig}
\efi

\ble[Three-Points Lemma]
\label{tpl}
Let $u<v<w$ be three points in the same equivalence class and let $P$
be a path.  Suppose that 
$P$ passes strictly west of $u$, weakly east of $v$, and weakly
west of $w$.
Then $P$ is invalid.
\ele
\pf\
Among all such triples $(u,v,w)$ satisfying the conditions of the
lemma, we can choose one where $u,v$ 
are a minimum distance apart, in which case
$v=u+(r,s)$.  Now from the possible $w$'s satisfying the hypotheses of
the lemma with this $u,v$, pick the one which has minimum distance
from $v$.   Let $i$ and $k$ be the integers such that $v$ and $w$
are on the lines $x+y=i+1$ and $x+y=k$, respectively.  Also let
$Q=P+(r,s)$.   Note that even though $P$ and $Q$ could intersect on
$x+y=i+1$, they could only do so if they entered with an $N$-step and
an $E$-step, respectively.  So we have $x_i(P)>x_i(Q)$.  It is also
clear from the choice of $w$ that $x_k(P)\le x_k(Q)$.  Thus we
are done by the Switching Lemma.
\Qqed

\bfi
\begin{pspicture}(0,0)(10,6)
\psgrid[subgriddiv=1,griddots=10](0,0)(10,6)
\psline(0,0)(0,1)
\psline(0,1)(4,1)
\psline(4,1)(4,6)
\psline(4,6)(9,6)
\psline[linestyle=dashed](0,0)(3,0)
\psline[linestyle=dashed](3,0)(3,2)
\psline[linestyle=dashed](3,2)(6,2)
\psline[linestyle=dashed](6,2)(6,4)
\psline[linestyle=dashed](6,4)(9,4)
\psline[linestyle=dashed](9,4)(9,6)
\pscircle*(1,0){.1}
\rput(1.3,.3){$u$}
\pscircle*(4,2){.1}
\rput(4.3,2.3){$v$}
\pscircle*(10,6){.1}
\rput(10.3,6.3){$w$}
\end{pspicture}
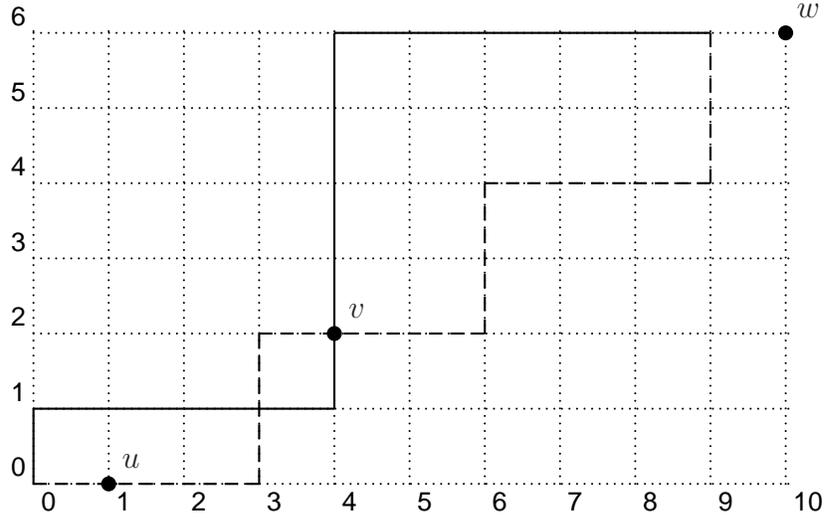
\capt{The Staircase Lemma for $r=3$, $s=2$, and $n=3$}
\label{stlfig}
\efi

Given two paths $P,Q$  we say that $P$ is
{\it northwest\/} of $Q$  if 
for every vertex $v=(x,y)$ of $P$ there is a vertex $w=(x',y')$ of
$Q$ southeast of $v$, i.e., $x\le x'$ and $y\ge y'$.  The {\it
staircase\/} is the path from $(0,0)$ to $(nr,ns)$ with steps
$$
S=(E^r N^s)^n.
$$
The dashed path in Figure~\ref{stlfig} is the staircase for $r=3$,
$s=2$, and $n=3$. 
The staircase forms a natural boundary for valid paths.  In following
the proof of the following lemma, the reader may wish to consult
Figure~\ref{stlfig}.  
\ble[Staircase Lemma]
If $P$ is a valid path from $(0,0)$ to $(nr,ns)$ then $P$ is northwest
of $S$.
\ele
\pf\
Suppose not. Then since $P$ ends northwest of $S$, we can pick a
vertex $v$ which is the first vertex of $P\cap S$ after some vertex of
$P$ which is (strictly) southeast of $S$.  Note that $P$ must enter
$v$ with an $N$-step and so $P$ passes weakly east of $v$.  Note also
that $S$ must enter $v$ with an $E$-step.  It
follows that we have points $u$ and $w$ in $[v]$ which are on the
lines $y=0$ and $y=ns$, respectively, but that these points are
strictly east of the corresponding points in $[0,0]$.  Since $P$ begins at $(0,0)$ and
ends at $(nr,ns)$ which are both in $[0,0]$, $P$ 
passes strictly west of $u$ and weakly west of $w$.  So $P$ is invalid
by the Three-Points Lemma, a contradiction.
\Qqed

\section{The bijection}

Let $\cV_n$ be the set of all
valid lattice paths from $(0,0)$ to $(nr,ns)$.  Then to prove
Conjecture~\ref{lw}, it suffices to find, for each $n\ge2$, a bijection
$\phi:\cV_n\ra\cV_1\times \cV_{n-1}$.  To see this, note that
every lattice path from $(0,0)$ to $(r,s)$ is valid and so
$|\cV_1|={r+s\choose r}$.  Iterating this map gives a bijection
between $\cV_n$ and $(\cV_1)^n$.  And the latter is clearly counted by
${r+s\choose r}^n$.

\bfi
\begin{pspicture}(0,0)(6,6)
\psgrid[subgriddiv=1,griddots=10](0,0)(6,6)
\psline(0,0)(0,1)
\psline(0,1)(1,1)
\psline(1,1)(1,5)
\psline(1,5)(4,5)
\psline(4,5)(4,6)
\psline(4,6)(6,6)
\pscircle*(1,2){.1}
\pscircle*(5,6){.1}
\rput(.5,.5){$P$}
\rput(1.5,4.5){$Q$}
\rput(3,2.5){$X$}
\rput(5,4.5){$X$}
\end{pspicture}
  \begin{pspicture}(0,0)(2,6)
  \rput(1,3){$\stackrel{\textstyle \phi}{\mapsto}$}
  \end{pspicture}
\begin{pspicture}(0,0)(6,6)
\psgrid[subgriddiv=1,griddots=10](0,0)(6,6)
\psline(0,0)(0,1)
\psline(0,1)(1,1)
\psline(1,1)(1,2)
\psline(1,2)(2,2)
\psline(2,2)(2,5)
\psline(2,5)(5,5)
\psline(5,5)(5,6)
\psline(5,6)(6,6)
\pscircle*(2,2){.1}
\pscircle*(6,6){.1}
\rput(.5,.5){$P'$}
\rput(2.5,4.5){$Q'$}
\end{pspicture}
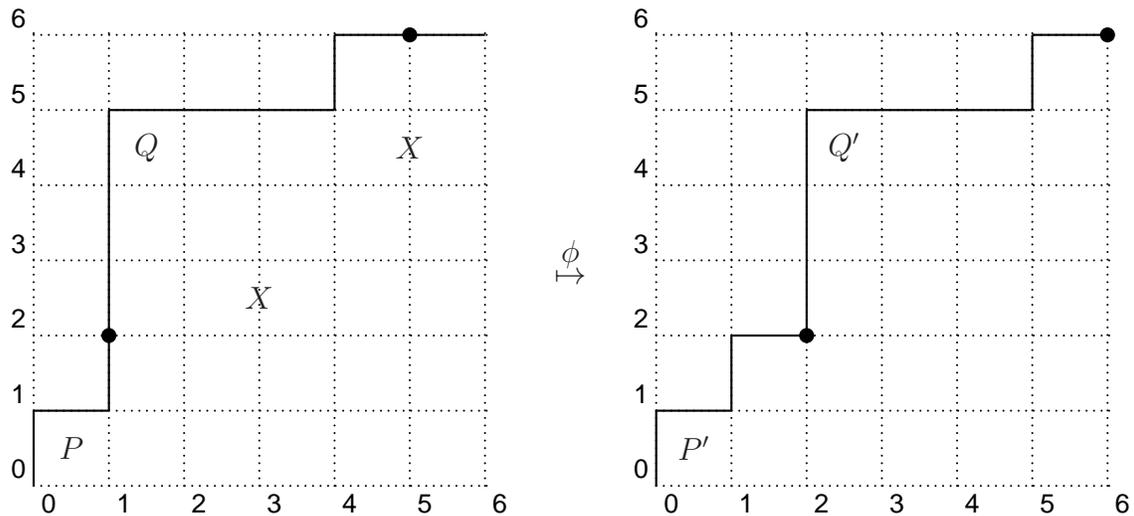
\capt{The first case of the bijection with $r=s=2$, $a=0$, and $b=c=1$}
\label{case1}
\efi

\noindent {\bf Proof (of Conjecture~\ref{lw} for $s=2$)}\quad
We construct the bijection $\phi$ when $s=2$.  Given $P\in\cV_n$ we wish
to construct $\phi(P)=(P',Q')\in\cV_1\times\cV_{n-1}$.  
By convention, we will consider $P'$ as going from $(0,0)$ to $(r,s)$
and $Q'$ as going from $(r,s)$ to $(nr,ns)$.  (Strictly speaking, $Q'$
is not in $\cV_{n-1}$ since it doesn't begin at the origin.  But the
translation of a valid path is valid, so no harm is done.)
Consider the prefix of $P$ up to and including the second $N$-step
which is of the form $E^a N E^b N$ for some $a,b\ge0$.  By the
Staircase Lemma, we must have $a+b\le r$.  Now consider 
the suffix of $P$ following the last $N$-step.  
Depending on whether the combined number of $E$-steps in the prefix
and the suffix is at least $r$ or less than $r$, we have two cases.
\ben
\item[(1)] For some path $Q$ we have
$$
P=E^a N E^b N Q E^c\qmq{where} a+b+c=r.
$$
\item[(2)] For some path $Q$ we have
$$
P=E^a N E^b N Q N E^c\qmq{where} a+b+c<r.
$$
\een

In the first case, we let
$$
P'=E^a N E^b N E^c\qmq{and} Q'=Q+(c,0).
$$
An example of such a path $P$ and its image is given in
Figure~\ref{case1} where the  circles indicate the endpoints
of $Q$ and $Q'$.
To show that the map and its inverse are well-defined, we will need the
following concept and result.  Given two paths $A$ and $B$, a {\it
blocked edge\/} is an $N$-edge or an $E$-edge of $\bbZ^2$ which can
not be on any valid path having $A$ as its prefix and $B$ as its
suffix.  Such edges will be marked with $X$'s in our figures.
\bpr
In case (1), the path $Q'$ is valid and begins with at most $r-b$ $E$-steps.
\epr
\pf\
Clearly $Q'$ is valid being a translate of a subpath of a valid path.

For the second statement, suppose to the contrary that $Q'$ begins
with more than $r-b$ $E$-steps.  Then $Q$ contains a point $(d,2)$
where $d>r+a$.  But now it is
impossible for $P$ to get to $(rn,2n)$.  This is because the initial
$E$-steps of $P$ produce a
sequence of blocked $N$-edges starting at points $(x,y(x))$ for every
$x$ with $d\le x\le rn$, where $2\le y(x)<sn$ and $y(x)$ is a weakly
increasing function of $x$.  So we have a contradiction.
\Qqed

\bfi
\begin{pspicture}(0,0)(6,6)
\psgrid[subgriddiv=1,griddots=10](0,0)(6,6)
\psline(0,0)(0,5)
\psline(0,5)(5,5)
\psline(5,5)(5,6)
\psline(5,6)(6,6)
\pscircle*(0,2){.1}
\pscircle*(5,5){.1}
\pscircle(0,2){.2}
\pscircle(2,5){.2}
\rput(.5,.5){$P$}
\rput(.5,4.5){$R$}
\rput(.5,2){$X$}
\rput(2.5,4){$X$}
\end{pspicture}
  \begin{pspicture}(0,0)(2,6)
  \rput(1,3){$\stackrel{\textstyle \phi}{\mapsto}$}
  \end{pspicture}
\begin{pspicture}(0,0)(6,6)
\psgrid[subgriddiv=1,griddots=10](0,0)(6,6)
\psline(0,0)(0,1)
\psline(0,1)(1,1)
\psline(1,1)(1,2)
\psline(1,2)(4,2)
\psline(4,2)(4,6)
\psline(4,6)(6,6)
\pscircle*(2,2){.1}
\pscircle*(6,6){.1}
\pscircle(4,3){.2}
\pscircle(6,6){.2}
\rput(.5,.5){$P'$}
\rput(4.5,5.5){$R'$}
\end{pspicture}
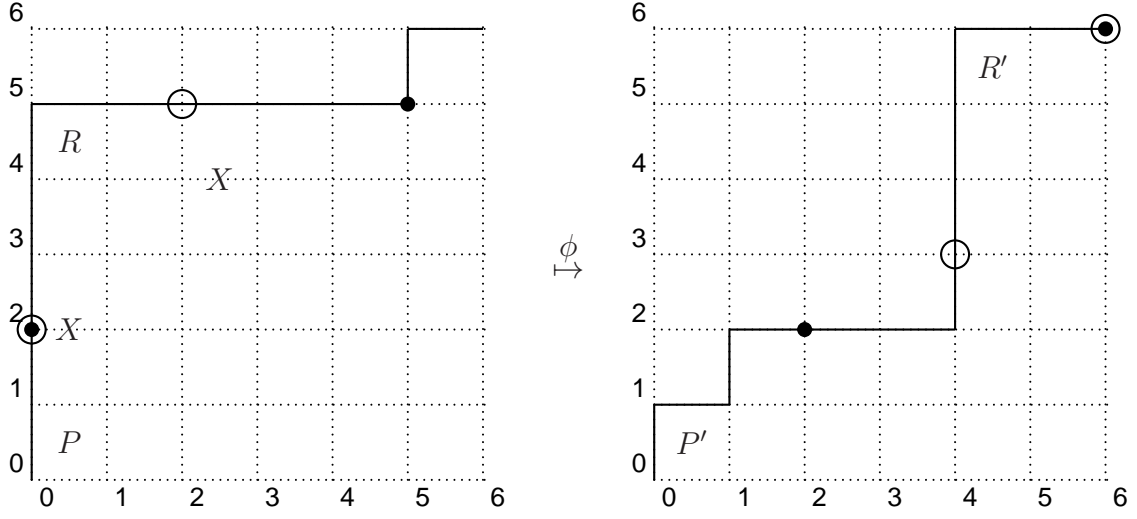
\capt{The second case of the bijection with $r=s=2$, $a=b=0$, and $c=1$}
\label{case2}
\efi

In the second case, we will show that $Q=R E^{r+1}$ for some path
$R$.  Assuming this for the moment, we can define
$$
P'=E^a N E^{r-a-c} N E^c\qmq{and} Q'=E^{a+b+c+1} N R'
$$
where $R'=R+(c+r+1,1)$.  Figure~\ref{case2} illustrates this case with
the endpoints of $Q$ and $Q'$ being marked with closed circles while
those of $R$ and $R'$ are marked with open ones.

We now verify the claim about $Q$ and the fact
that $Q'$ is valid.
\bpr
In case (2), the path $Q$ ends with at least $r+1$ $E$-steps.  In addition, $Q'$ is valid.
\epr 
\pf\
We prove the first statement by contradiction.  Note that because of
the final $N$-step in $P$, any $E$-edge into a vertex of the class
$[r-c,2]$ is blocked.  Now  the subpath $QN$  of $P$ passes strictly
west of $u=(r-c,2)$ since $a+b+c<r$.  This subpath also ends at $w=u+(n-1)(r,2)$.
But if $Q$ ends 
with fewer than $r+1$ $E$-steps then $QN$ passes weakly east of 
$v=u+(n-2)(r,2)$ since the $E$-edge into $v$ is
blocked.  This contradicts the Three-Points Lemma as long as 
$n\ge3$.  For $n=2$, just note that $P$ would have to contain points
both east and west of the blocked edge into $(r-c,2)$ which is impossible.

Since $R'$ is a translation of a valid path $R$, it is valid itself.
So the only way $Q'$ could be invalid is if one of the $E$-steps in
the prefix $E^{a+b+c+1}$ is in conflict with an $N$-step in $R'$.
Note that such an $N$-step must be out of a point of some class $[x,3]$
with $x\ge 2r+1$.
Suppose that this is the case and consider what this implies about the
original path $P$.  In particular, consider the class
$[r-c,2]$ as in the previous paragraph.  Using the same $u$ and $w$ as
before, the supposed 
$N$-step forces $P$ to contain a point weakly east
of some $v$ in this class with
$u<v<w$.  But then $P$ must pass weakly east of $v$ since the $E$-edge
into $v$ is blocked.
So the Three-Points Lemma (or a direct argument when $n=2$)
provides the necessary contradiction. 
\Qqed

We now describe the inverse map.   Suppose we are given
$(P',Q')\in\cV_1\times\cV_{n-1}$ and write
$$
P'=E^a N E^b N E^c\qmq{and} Q'=E^dNR'
$$
for some path $R'$.  Then, again, we have two cases to describe
$P=\phi^{-1}(P',Q')$. 
\ben
\item  If $d\le r-b$ then let
$$
P=E^a N E^b N Q E^c\qmq{where} Q = Q'-(c,0).
$$
\item  If $d> r-b$ then let
$$
P=E^a N E^{b+d-r-1} N R E^{r+1} N E^c\qmq{where} R=R'-(r+c+1,1).
$$
\een

It is easy to verify that this is a case-by-case inverse for the map
$\phi$.  Furthermore, the demonstration that $\phi^{-1}$ is well
defined is quite similar to the one just given for $\phi$, so we omit
it.  This completes the proof of Conjecture~\ref{lw} when $s=2$.
\Qqed

We can say a little more about the case $s=2$.  Let
$\Phi:\cV_n\ra(\cV_1)^n$ be the map obtained by composing $\phi$ with
itself $n-1$ times.  Consider a path $P\in\cV_n$ and let
$\Phi(P)=(P_1',\ldots,P_n')$.  Then directly from our definition of
$\phi$, we see that $P$ and $P_1'$ begin with the same number of
$E$-steps before the first $N$-step.  So given $a$ with $0\le a\le r$, 
$\Phi$ restricts to a bijection between the $P$ with prefix $E^a N$
and the $n$-tuples $(P_1',\ldots,P_n')$ where
$P_1'$ satisfies the same restriction.  But, as we mentioned before,
the validity condition imposes no restriction on paths in $\cV_1'$, so
the number of such $P_1'$ is clearly $r-a+1$.  Thus we have proved the
following corollary.
\bco
Suppose $s=2$.  Given $a$ with $0\le a\le r$, the number of
$P\in\cV_n$ with a prefix of the form $E^aN$ is 
$(r-a+1){r+2\choose 2}^{n-1}$.\Qqed
\eco

\bigskip
\bibliographystyle{acm}
\begin{small}
\bibliography{ref}
\end{small}

\end{document}